\documentclass[twoside,12pt,
a4paper,backrefs,msc-links]{amsart}
\usepackage{graphicx}
\PassOptionsToPackage{pdfauthor={Vladimir V. Kisil},%
    pdftitle={Erlangen Program at Large---0: Starting with the group SL(2,R)},%
    pdfsubject={general mathematics},%
    backref=page,%
    breaklinks=true,
    pdfkeywords={symmetries, Erlangen program}
}{hyperref}
\usepackage{hyperref}
\usepackage{multicol}
\usepackage{amsrefs}
\IfFileExists{eulervm.sty}{\usepackage{eulervm}
}{}
\newtheorem{thm}{Theorem}
\theoremstyle{definition}
\newtheorem{defn}[thm]{Definition}

\providecommand{\lvec}[1]{\overrightarrow{#1}}
\providecommand{\cycle}[3][]{{#1 C^{#2}_{#3}}}
\newcommand{\zcycle}[3][]{#1 Z^{#2}_{#3}}
\newcommand{\realline}[3][]{#1 R^{#2}_{#3}}
\providecommand{\GiNaC}{\textsf{GiNaC}}
\providecommand{\bs}{\breve{\sigma}}
\providecommand{\SL}[1][2]{\ensuremath{\FSpace{SL}{#1}(\Space{R}{})}}
\providecommand{\scalar}[3][\relax]{\left\langle #2,#3 
        \right\rangle\ifx#1\relax\else_{#1}\fi}
\providecommand{\Space}[3][]{\ensuremath{\mathbb{#2}^{#3}_{#1}{}}}
  \providecommand{\FSpace}[3][]{\ensuremath{\ifx#2l \ell_{#3}^{#1}{}\else
  #2_{#3}^{#1}{}\fi}} 
\providecommand{\rmi}{\mathrm{i}}
\providecommand{\rmc}{\mathrm{\breve\i}}
\providecommand{\tr}{\mathop{tr}}
\providecommand{\algebra}[1]{\ensuremath{\mathfrak{#1}}}
\providecommand{\MR}[1]{\textbf{MR}~\href{http://www.ams.org/mathscinet-getitem?mr=#1}{\#~#1}}

\providecommand{\eprint}[2]{E-print: \href{#1}{\texttt{#2}}}
\providecommand{\modulus}[2][\relax]{\left| #2 \right|\ifx#1\relax\else_{#1}\fi}
\hypersetup{colorlinks=true,bookmarks=true}
\providecommand{\wiki}[2]{\href{http://en.wikipedia.org/wiki/#1}{#2}}

\begin{document}
\title[EPAL0: Starting with SL(2,R) group]{Erlangen Program at
  Large---0: Starting with the Group \(\SL\) }

\author[Vladimir V. Kisil]%
{\href{http://maths.leeds.ac.uk/~kisilv/}{Vladimir V. Kisil}}
\thanks{On  leave from the Odessa University.}
\dedicatory{Dedicated to the memory of Serge Lang}
\address{School of Mathematics\\
University of Leeds\\
Leeds LS2\,9JT\\
UK
}

\email{\href{mailto:kisilv@maths.leeds.ac.uk}{kisilv@maths.leeds.ac.uk}}

\urladdr{\href{http://maths.leeds.ac.uk/~kisilv/}%
{http://maths.leeds.ac.uk/\~{}kisilv/}}

\maketitle

\begin{multicols}{2}
  
  The simplest objects with non-commutative multiplication may be
  \(2\times 2\) matrices with real entries.  Such matrices \emph{of
    determinant one} form a closed set under multiplication (since
  \(\det (AB)=\det A\cdot \det B\)), the identity matrix is among them
  and any such matrix has an inverse (since \(\det A\neq 0\)). In
  other words those matrices form a group, \wiki{SL2(R)}{the \(\SL\)
    group}~\cite{Lang85}---one of the two most important Lie groups
  in analysis. The other group is \wiki{Heisenberg_group}{the
    Heisenberg group}~\cite{Howe80a}. By contrast the
  \wiki{Affine_transformation}{``\(ax+b\)''-group}, which is often
  used to build wavelets, is only a subgroup of \(\SL\), see the
  numerator in~\eqref{eq:moebius}.

  The simplest non-linear transforms of the real
  line---linear-fractional or \wiki{Moebius_transformation}{M\"obius
    maps}---may also be associated with \(2\times 2\)
  matrices~\cite{Beardon05a}*{Ch.~13}:
  \begin{equation}
    \label{eq:moebius}
    g: x\mapsto g\cdot x=\frac{ax+b}{cx+d}, \text{ where } 
    g=  \begin{pmatrix}
      a&b\\c&d
    \end{pmatrix}, x\in\Space{R}{}.
  \end{equation}
  An enjoyable calculation shows that the composition of two
  transforms~\eqref{eq:moebius} with different matrices \(g_1\) and
  \(g_2\) is again a M\"obius transform with matrix the product
  \(g_1 g_2\). In other words~\eqref{eq:moebius} it is a (left) action
  of \(\SL\).

  According to F.~Klein's \wiki{Erlangen_program}{\emph{Erlangen
      program}} (which was influenced by S.~Lie) any geometry is
  dealing with invariant properties under a certain group action. For
  example, we may ask: \emph{What kinds of geometry are related to
   the \(\SL\) action~\eqref{eq:moebius}}?
  
  The Erlangen program has probably the highest rate of
  \(\frac{\text{praised}}{\text{actually used}}\) among mathematical
  theories not only due to the big numerator but also due to undeserving
  small denominator. As we shall see below Klein's approach provides
  some surprising conclusions even for such over-studied objects as
  circles.

  \section{Make a Guess in Three Attempts}
  \label{sec:make-guess-three}

  It is easy to see that the \(\SL\) action~\eqref{eq:moebius} makes
  sense also as a map of complex numbers \(z=x+\rmi y\),
  \(\rmi^2=-1\). Moreover, if \(y>0\) then \(g\cdot z\) has a positive
  imaginary part as well, i.e. \eqref{eq:moebius} defines a map from
  the upper half-plane to itself. 

  However there is no need to be restricted to the traditional route
  of complex numbers only. Less-known \wiki{Dual_number}{\emph{dual}}
  and \wiki{Split-complex_number}{\emph{double}} numbers
  \cite{Yaglom79}*{Suppl.~C} have also the form \(z=x+\rmi y\) but
  different assumptions on the imaginary unit \(\rmi\): \(\rmi^2=0\)
  or \(\rmi^2=1\) correspondingly. Although the arithmetic of dual and
  double numbers is different from the complex ones, e.g. they have
  divisors of zero, we are still able to define their transforms
  by~\eqref{eq:moebius} in most cases.

  Three possible values \(-1\), \(0\) and \(1\) of \(\sigma:=\rmi^2\)
  will be refereed to here as \emph{elliptic}, \emph{parabolic} and
  \emph{hyperbolic} cases respectively.  We repeatedly meet such a
  division of various mathematical objects into three classes.  They
  are named by the historically first example---the classification of
  conic sections---however the pattern persistently reproduces itself
  in many different areas: equations, quadratic forms, metrics,
  manifolds, operators, etc.  We will abbreviate this separation as
  \emph{EPH-classification}.  The \emph{common origin} of this
  fundamental division can be seen from the simple picture of a
  coordinate line split by zero into negative and positive
  half-axes:
  \begin{equation}
    \label{eq:eph-class}
    \raisebox{-15pt}{\includegraphics[scale=1]{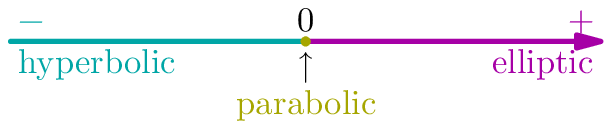}}
  \end{equation}

  Connections between different objects admitting EPH-classification
  are not limited to this common source. There are many deep results
  linking, for example,
  \wiki{Atiyah-Singer_index_theorem}{ellipticity
    of quadratic forms, metrics and operators}.
 On the other hand there are still a lot of white spots and obscure
 gaps between some subjects as well.

  To understand the action~\eqref{eq:moebius} in all EPH cases we use
  the Iwasawa decomposition~\cite{Lang85} of \(\SL=ANK\) into
  \emph{three} one-dimensional subgroups \(A\), \(N\) and \(K\):
  \begin{equation}
    \label{eq:iwasawa-decomp}
    \begin{pmatrix}
      a&b \\c &d
    \end{pmatrix}= {\begin{pmatrix} \alpha & 0\\0&\alpha^{-1}
      \end{pmatrix}} {\begin{pmatrix} 1&\nu \\0&1
      \end{pmatrix}} {\begin{pmatrix}
        \cos\phi &  \sin\phi\\
        -\sin\phi & \cos\phi
      \end{pmatrix}}.
  \end{equation}
  Subgroups \(A\) and \(N\) act in~\eqref{eq:moebius} irrespectively
  to value of \(\sigma\): \(A\) makes a dilation by \(\alpha^2\), i.e.
  \(z\mapsto \alpha^2z\), and \(N\) shifts points to left by \(\nu\),
  i.e. \(z\mapsto z+\nu\).

\begin{figure}[htbp]
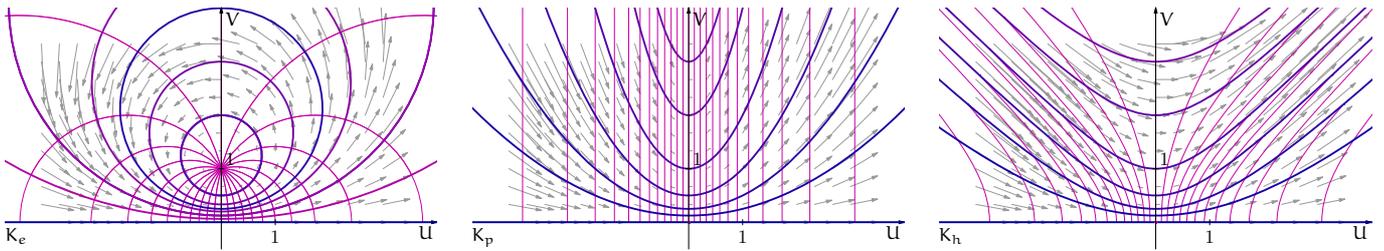

  \centering
  \includegraphics[scale=.71]{parabolic0.4}\hfill
  \includegraphics[scale=.71]{parabolic0.5}\hfill
  \includegraphics[scale=.71]{parabolic0.6}
  \caption[Action of the $K$ subgroup]{Action of the \(K\) subgroup.
    The corresponding \(K\)-orbits are thick circles, parabolas and
    hyperbolas. Thin traversal lines are images of the vertical axis
    for certain values of the parameter \(\phi\).}
  \label{fig:k-subgroup}
\end{figure}

  By contrast, the action of the third matrix from the subgroup \(K\)
  sharply depends on \(\sigma\), see Fig.~\ref{fig:k-subgroup}. In
  elliptic, parabolic and hyperbolic cases \(K\)-orbits are circles,
  parabolas and (equilateral) hyperbolas correspondingly.  Thin
  traversal lines in Fig.~\ref{fig:k-subgroup} join points of orbits
  for the same values of \(\phi\) and grey arrows represent ``local
  velocities''---vector fields of derived representations.

  \begin{defn}
    \label{de:cycle}
    The common name \emph{cycle}~\cite{Yaglom79} is used to
    denote circles, parabolas and hyperbolas (as well as straight lines
    as their limits) in the respective EPH case.
  \end{defn}

\begin{figure}[htbp]
  \centering
  (a)\includegraphics{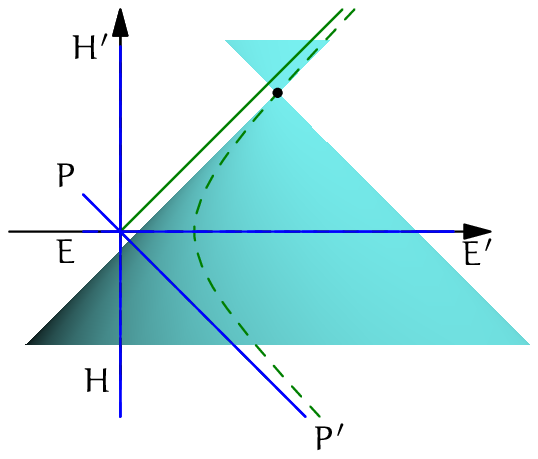}\hspace{2cm}
  (b)\includegraphics{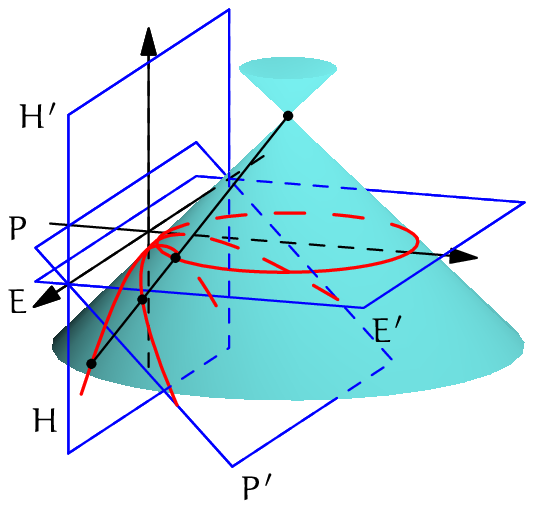}
  \caption[$K$-orbits as conic sections]{\(K\)-orbits as conic
    sections: 
    circles are sections by the plane \(EE'\);  parabolas are
    sections by \(PP'\);  hyperbolas are sections by \(HH'\). Points
    on the same generator of the cone correspond to the same value of \(\phi\).}
  \label{fig:k-orbit-sect}
\end{figure}

\begin{figure}[htpb]
  \centering
  \includegraphics[scale=.75]{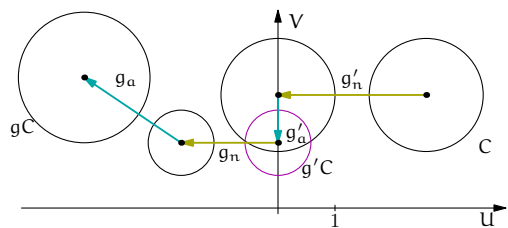}
  \caption[The decomposition of an arbitrary Moebius
  transformation]{Decomposition of an arbitrary M\"obius
    transformation \(g\) into a product \(g=g_a g_n g_k g_a' g_n'\).}
  \label{fig:moeb-decomp}
\end{figure}

  It is well known that any cycle is a \emph{conic sections} and an
  interesting observation is that corresponding \(K\)-orbits are in
  fact sections of the same two-sided right-angle cone, see
  Fig.~\ref{fig:k-orbit-sect}.  Moreover, each straight line
  generating the cone, see Fig.~\ref{fig:k-orbit-sect}(b), is crossing
  corresponding EPH \(K\)-orbits at 
  points with the same value of parameter \(\phi\)
  from~\eqref{eq:iwasawa-decomp}. In other words, all three types of
  orbits are generated by the rotations of this generator along the
  cone.

  \(K\)-orbits are \(K\)-invariant in a trivial way. Moreover since
  actions of  both \(A\) and \(N\) for any \(\sigma\) are extremely
  ``shape-preserving'' we find natural invariant objects of the
  M\"obius map: 

  \begin{thm}
    The family of all cycles from Defn.~\ref{de:cycle} is invariant under the
    action~\eqref{eq:moebius}.
  \end{thm}
  \begin{proof}
    We will show that for a given \(g\in \SL\) and a cycle \(C\) its
    image \(gC\) is again a cycle. Fig.~\ref{fig:moeb-decomp} make an
    illustration with \(C\) as a circle, but our reasoning works in
    all EPH cases.

    For a fixed \(C\) there is always the unique pair of
    transformations \(g'_n\) from the subgroup \(N\) and \(g'_a\in A\)
    that the cycle \(g'_a g'_n C\) is exactly a \(K\)-orbit. We make a
    decomposition of \(g (g'_a g'_n)^{-1}\) into a product
    as in~\eqref{eq:iwasawa-decomp}:
    \begin{displaymath}
      g (g'_a g'_n)^{-1} = g_ag_ng_k.
    \end{displaymath}
    Since \(g'_a g'_n
    C\) is a \(K\)-orbit we have \(g_k (g'_a g'_n C)=g'_a g'_n C\),
    then:
    \begin{eqnarray*}
      gC&=& g (g'_a g'_n)^{-1} g'_a g'_n C = g_ag_n g_k g'_a g'_n C \\
      &=& g_ag_n g_k(g'_a g'_n C) =  g_ag_n g'_a g'_n C,
    \end{eqnarray*}
    Since subgroups \(A\) and \(N\) obviously preserve the shape of
    any cycle  this finishes our proof.
  \end{proof}
  According to Erlangen ideology we should now study invariant
  properties of cycles. 

  \section{Invariance of FSCc}
  \label{sec:invariance-fscc}

  Fig.~\ref{fig:k-orbit-sect} suggests that we may get a unified
  treatment of cycles in all EPH by consideration of a higher
  dimension spaces. The standard mathematical method is to declare
  objects under investigations (cycles in our case, functions in
  functional analysis, etc.) to be simply points of some bigger
  space. This space should be equipped with an appropriate structure
  to hold externally information which were previously inner
  properties of our objects.
  
  A generic cycle is the set of points \((u,v)\in\Space{R}{2}\)
  defined for all values of \(\sigma\) by the equation
  \begin{equation}
    \label{eq:cycle-eq}
    k(u^2-\sigma v^2)-2lu-2nv+m=0.
  \end{equation}
  This equation (and the corresponding cycle) is defined by a point
  \((k, l, n, m)\) from a projective space \(\Space{P}{3}\), since for
  a scaling factor \(\lambda \neq 0\) the point \((\lambda k, \lambda
  l, \lambda n, \lambda m)\) defines the same
  equation~\eqref{eq:cycle-eq}. We call \(\Space{P}{3}\) the
  \emph{cycle space} and refer to the initial \(\Space{R}{2}\) as the
  \emph{point space}.

  In order to get a connection with M\"obius action~\eqref{eq:moebius}
  we arrange numbers \((k, l, n, m)\) into the matrix 
  \begin{equation}
    \label{eq:FSCc-matrix}
    C_{\bs}^s=\begin{pmatrix}
      l+\rmc s n&-m\\k&-l+\rmc s n
    \end{pmatrix}, 
  \end{equation}
  with a new imaginary unit \(\rmc\) and an additional parameter \(s\)
  usually equal to \(\pm 1\). The values of \(\bs:=\rmc^2\) is \(-1\),
  \(0\) or \(1\) independently from the value of \(\sigma\).  The
  matrix~\eqref{eq:FSCc-matrix} is the cornerstone of (extended)
  Fillmore--Springer--Cnops construction (FSCc)~\cite{Cnops02a} and
  closely related to technique recently used by A.A.~Kirillov to study
  the Apollonian gasket~\cite{Kirillov06}.

  The significance of FSCc in Erlangen framework is provided by the
  following result:
  \begin{thm}
    \label{th:FSCc-intertwine}
    The image  \(\tilde{C}_{\bs}^s\) of a cycle \(C_{\bs}^s\) under
    transformation~\eqref{eq:moebius} with \(g\in\SL\) is given by
    similarity of the matrix~\eqref{eq:FSCc-matrix}:
    \begin{equation}
      \label{eq:cycle-similarity}
      \tilde{C}_{\bs}^s= gC_{\bs}^sg^{-1}.
    \end{equation}
    In other words FSCc~\eqref{eq:FSCc-matrix} \emph{intertwines}
    M\"obius action~\eqref{eq:moebius} on cycles with
    linear map~\eqref{eq:cycle-similarity}.
  \end{thm}
  
  There are several ways to prove~\eqref{eq:cycle-similarity}: either
  by a brute force calculation (fortunately
  \href{http://arxiv.org/abs/cs.MS/0512073}{performed by a
    CAS})~\cite{Kisil05a} or through the related orthogonality of
  cycles~\cite{Cnops02a}, see the end of the next
  section~\ref{sec:invar-algebr-geom}. 

  The important observation here is that FSCc~\eqref{eq:FSCc-matrix}
  uses an imaginary unit \(\rmc\) which is not related to \(\rmi\)
  defining the appearance of cycles on plane. In other words any EPH
  type of geometry in the cycle space \(\Space{P}{3}\) admits drawing
  of cycles in the point space \(\Space{R}{2}\) as circles, parabolas
  or hyperbolas. We may think on points of \(\Space{P}{3}\) as ideal
  cycles while their depictions on \(\Space{R}{2}\) are only their
  shadows on the wall of
  \href{http://en.wikipedia.org/wiki/Plato#Metaphysics}{Plato's
    cave}. 

\begin{figure}[htbp]
  \centering
  (a) \includegraphics[]{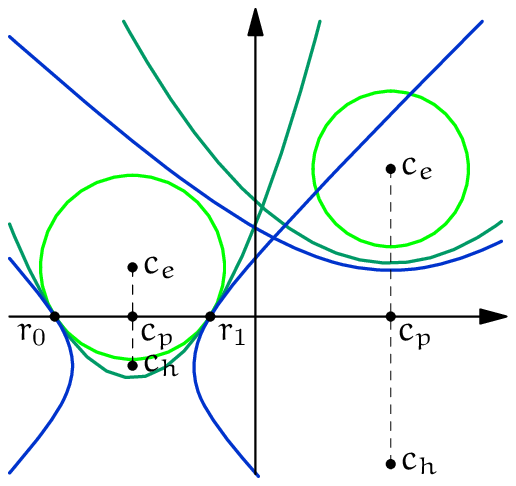}\hspace{1cm}
  (b) \includegraphics[]{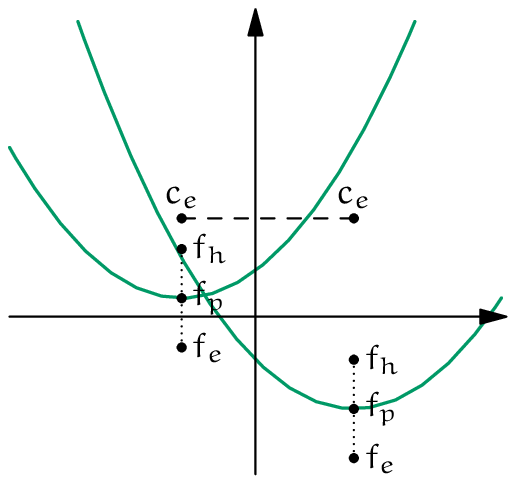}
  \caption[Cycle implementations, centres and foci]{
    (a) Different
    EPH implementations of the same cycles defined by quadruples of
    numbers.\\
    (b) Centres and foci of two parabolas with the same focal length.} 
  \label{fig:eph-cycle}
\end{figure}

  Fig.~\ref{fig:eph-cycle}(a) shows the same cycles drawn
  in different EPH styles. Points \(c_{e,p,h}=(\frac{l}{k}, -\sigma
  \frac{n}{k})\) are their respective e/p/h-centres. They are related to
  each other through several identities:
  \begin{equation}
    \label{eq:centres}
    c_e=\bar{c}_h, \quad c_p=\frac{1}{2}(c_e+c_h).
  \end{equation}
  Fig.~\ref{fig:eph-cycle}(b) presents two cycles drawn as parabolas,
  they have the same focal length \(\frac{n}{2k}\) and thus their
  e-centres are on the same level. In other words \emph{concentric}
  parabolas are obtained by a vertical shift, not
  scaling as an analogy with circles or hyperbolas may suggest. 

  Fig.~\ref{fig:eph-cycle}(b) also presents points, called e/p/h-foci:
  \begin{equation}
    \label{eq:foci}
    f_{e,p,h}=\left(\frac{l}{k}, -\frac{\det C_{\bs}^s}{2nk}\right),
  \end{equation}
  which are independent of the sign of \(s\).  If a cycle is depicted
  as a parabola then h-focus, p-focus, e-focus are correspondingly
  geometrical focus of the parabola, its vertex, and the point on the
  directrix nearest to the vertex.

  As we will see, cf. Thms.~\ref{th:ghost1} and~\ref{th:ghost2}, all
  three centres and three foci are useful  
  attributes of a cycle even if it is drawn as a circle.

  \section{Invariants: algebraic and geometric}
  \label{sec:invar-algebr-geom}

  We use known algebraic invariants of matrices to build appropriate
  geometric invariants of cycles. It is yet another demonstration that
  any division of mathematics into subjects is only illusive.

  For \(2\times 2\) matrices (and thus cycles) there are only two
  essentially different invariants under
  similarity~\eqref{eq:cycle-similarity} (and thus under M\"obius
  action~\eqref{eq:moebius}): the \emph{trace} and the
  \emph{determinant}.  The latter was already used in~\eqref{eq:foci}
  to define cycle's foci. However due to projective nature of the
  cycle space \(\Space{P}{3}\) the absolute values of trace or
  determinant are irrelevant, unless they are zero.

  Alternatively we may have a special arrangement for normalisation of
  quadruples \((k,l,n,m)\). For example, if \(k\neq0\) we may
  normalise the quadruple to
  \((1,\frac{l}{k},\frac{n}{k},\frac{m}{k})\) with highlighted cycle's
  centre. Moreover in this case \(\det \cycle{s}{\bs}\) is equal to
  the square of cycle's radius, cf. Section~\ref{sec:dist-lenght-perp}.
  Another normalisation \(\det \cycle{s}{\bs}=1\) is used
  in~\cite{Kirillov06} to get a nice condition for touching circles.

  We still get important characterisation even with non-normalised
  cycles, e.g., invariant classes (for different \(\bs\)) of
  cycles are defined by the condition \(\det C_{\bs}^s=0\). Such a
  class is parametrises only by two real number and as such is easily
  attached to certain point of \(\Space{R}{2}\). For example, the
  cycle \(C_{\bs}^s\) with \(\det C_{\bs}^s=0\), \(\bs=-1\) drawn
  elliptically represent just a point \((\frac{l}{k},\frac{n}{k})\),
  i.e. (elliptic) zero-radius circle.  The same condition with
  \(\bs=1\) in hyperbolic drawing produces a null-cone originated at
  point \((\frac{l}{k},\frac{n}{k})\):
  \begin{displaymath}
    (u-\frac{l}{k})^2-(v-\frac{n}{k})^2=0,
  \end{displaymath}
  i.e. a zero-radius cycle in hyperbolic metric. 

\begin{figure}[htbp]
  \centering
  \includegraphics{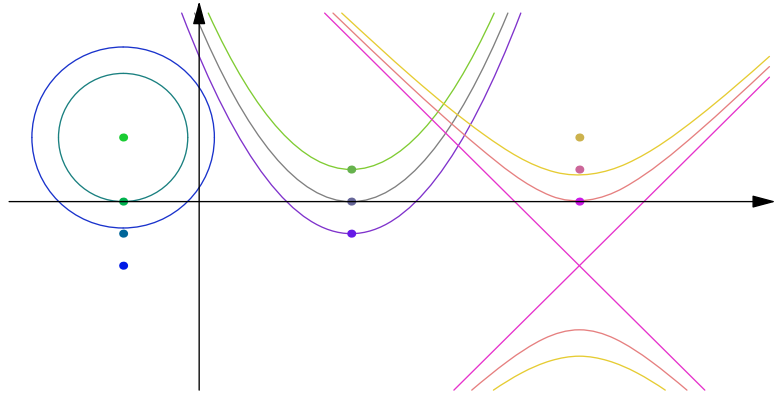}
  \caption[Different implementations of the same
    zero-radius cycles]{Different \(\rmi\)-implementations of the same
    \(\bs\)-zero-radius cycles and corresponding foci.}
  \label{fig:zero-radius}
\end{figure}

  In general for every notion there is nine possibilities: three EPH
  cases in the cycle space times three EPH realisations in the point
  space. Such nine cases for ``zero radius'' cycles is shown on
  Fig.~\ref{fig:zero-radius}. For example, p-zero-radius cycles in any
  implementation touch the real axis.

  This ``touching'' property is a manifestation of the \emph{boundary
    effect} in the upper-half plane
  geometry~\cite{Kisil05a}*{Rem.~3.4}. The famous question
  \href{http://en.wikipedia.org/wiki/Hearing_the_shape_of_a_drum}{on
    hearing drum's shape} has a sister:
  \begin{quote}
  \emph{Can we see/feel the boundary from inside a domain?}
  \end{quote}
  Both orthogonality relations  described below are ``boundary
  aware'' as well. It is not surprising after all since \(\SL\) action
  on the upper-half plane was obtained as an extension of its
  action~\eqref{eq:moebius} on the boundary. 

  According to the \wiki{Category_theory}{categorical viewpoint}
  internal properties of objects are of minor importance in comparison
  to their relations with other objects from the same class. As an
  illustration we may put the proof of Thm.~\ref{th:FSCc-intertwine}
  sketched at the end of of the next section. Thus from now on we will
  look for invariant relations between two or more cycles.

  \section{Joint invariants: orthogonality}
  \label{sec:joint-invar-orth}

  The most expected relation between cycles is based on the following
  M\"obius invariant ``inner product'' build from a trace of
  product of two cycles as matrices:
  \begin{equation}
    \label{eq:inner-prod}
    \scalar{C_{\bs}^s}{\tilde{C}_{\bs}^s}= \tr (C_{\bs}^s\tilde{C}_{\bs}^s)
  \end{equation}
  By the way, an inner product of this type is used, for example, in
  \wiki{Gelfand-Naimark-Segal_construction}{GNS construction} to make
  a Hilbert space out of \(C^*\)-algebra.  The next standard move is
  given by the following definition.
  \begin{defn}
    \label{de:orthogonality}
    Two cycles are called \(\bs\)-orthogonal if
    \(\scalar{C_{\bs}^s}{\tilde{C}_{\bs}^s}=0\). 
  \end{defn}
  For the case of \(\bs \sigma=1\), i.e. when geometries of the cycle
  and point spaces are both either elliptic or hyperbolic, such an
  orthogonality is the standard one, defined in terms of angles
  between tangent lines in the intersection points of two cycles.
  However in the remaining seven (\(=9-2\)) cases the innocent-looking
  Defn.~\ref{de:orthogonality} brings unexpected relations.

\begin{figure}[htbp]
  \includegraphics[]{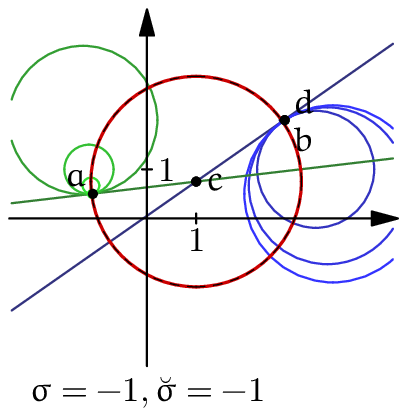}\hspace{1cm}
  \includegraphics[]{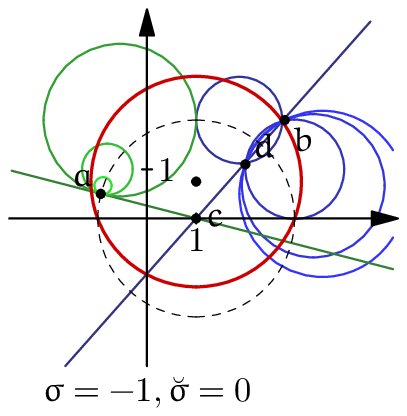}\hspace{1cm}
  \includegraphics[]{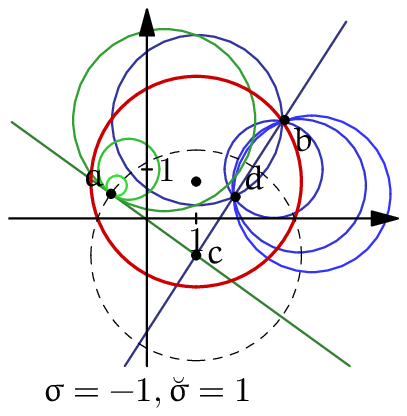}  \caption[Orthogonality of the first kind]{Orthogonality of the first
    kind in the elliptic point space.\\
    Each picture presents two groups (green and blue) of cycles which
    are orthogonal to the red cycle \(C^{s}_{\bs}\).  Point \(b\)
    belongs to \(C^{s}_{\bs}\) and the family of blue cycles
    passing through \(b\) is orthogonal to \(C^{s}_{\bs}\). They
    all also intersect in the point \(d\) which is the inverse of
    \(b\) in \(C^{s}_{\bs}\). Any orthogonality is reduced to the usual
    orthogonality with a new (``ghost'') cycle (shown by the dashed
    line), which may or may not coincide with \(C^{s}_{\bs}\). For
    any point \(a\) on the ``ghost'' cycle the orthogonality is
    reduced to the local notion in the terms of tangent lines at the
    intersection point. Consequently such a point \(a\) is always the
    inverse of itself.}
  \label{fig:orthogonality1}
\end{figure}

  Elliptic (in the point space) realisations of
  Defn.~\ref{de:orthogonality}, i.e. \(\sigma=-1\) is shown in
  Fig.~\ref{fig:orthogonality1}. The left picture corresponds to the
  elliptic cycle space, e.g. \(\bs=-1\). The orthogonality between
  the red circle and any circle from the blue or green families is
  given in the usual Euclidean sense. The central (parabolic in the
  cycle space) and the right (hyperbolic) pictures show non-local
  nature of the orthogonality.  There are analogues pictures in
  parabolic and hyperbolic point spaces as well~\cite{Kisil05a}.

  This orthogonality may still be expressed in the traditional sense
  if we will associate to the red circle the corresponding ``ghost''
  circle, which shown by the dashed line in Fig.~\ref{fig:orthogonality1}.
  To describe ghost cycle we need the
  \wiki{Heaviside_step_function}{\emph{Heaviside function}} \(\chi(\sigma)\):
  \begin{equation}
    \label{eq:heaviside-function}
    \chi(t)=\left\{
      \begin{array}{ll}
        1,& t\geq 0;\\
        -1,& t<0.
      \end{array}\right.
  \end{equation}

  \begin{thm}
    \label{th:ghost1}
    A cycle is \(\bs\)-orthogonal to cycle \(C_{\bs}^s\) if it is
    orthogonal in the usual sense to the \(\sigma\)-realisation of
    ``ghost'' cycle \(\hat{C}_{\bs}^s\), which is defined by the
    following two conditions:
    \begin{enumerate}
    \item \label{item:centre-centre-rel}
      \(\chi(\sigma)\)-centre of \(\hat{C}_{\bs}^s\) coincides
      with  \(\bs\)-centre of \(C_{\bs}^s\).
    \item Cycles \(\hat{C}_{\bs}^s\) and \(C^{s}_{\bs}\) have the same
      roots, moreover \(\det \hat{C}_{\sigma}^1= \det C^{\chi(\bs)}_{\sigma}\).
    \end{enumerate}
  \end{thm}
  The above connection between various centres of cycles illustrates
  their meaningfulness within our approach.

  One can easy check the following orthogonality properties of the
  zero-radius cycles defined in the previous section:
  \begin{enumerate}
  \item Since \(\scalar{C_{\bs}^s}{{C}_{\bs}^s}=\det {C}_{\bs}^s\)
    zero-radius cycles are  self-orthogonal (isotropic) ones.
  \item \label{it:ortho-incidence}
    A cycle \(\cycle{s}{\bs}\) is \(\sigma\)-orthogonal to a zero-radius
    cycle \(\zcycle{s}{\bs}\) if and only if \(\cycle{s}{\bs}\) passes
    through the \(\sigma\)-centre of \(\zcycle{s}{\bs}\).
  \end{enumerate}
  \begin{proof}[Sketch of proof of Thm.~\ref{th:FSCc-intertwine}]
    The validity of Thm.~\ref{th:FSCc-intertwine} for a zero-radius
    cycle 
    \begin{displaymath}
      \zcycle{s}{\bs}=
      \begin{pmatrix}
        z&-z\bar{z}\\1&-\bar{z}
      \end{pmatrix}=    
      \frac{1}{2}
      \begin{pmatrix}
        z&z\\1&1
      \end{pmatrix}    
      \begin{pmatrix}
      1&-\bar{z}\\1&-\bar{z}
    \end{pmatrix}      
  \end{displaymath}
  with the centre \(z=x+\rmi y\) is straightforward. This implies the
  result for a generic cycle with the help of M\"obius invariance of
  the product~\eqref{eq:inner-prod} (and thus the orthogonality) and
  the above relation~(\ref{it:ortho-incidence}) between the
  orthogonality and the incidence. See~\cite{Cnops02a} for details.
  \end{proof}

  \section{Higher order joint invariants: s-orthogonality}
  \label{sec:higher-order-joint}

  With appetite already wet one may wish to build more joint
  invariants. Indeed for any homogeneous polynomial
  \(p(x_1,x_2,\ldots,x_n)\) of several non-commuting variables one may
  define an invariant joint disposition of \(n\) cycles
  \({}^j\!\cycle{s}{\bs}\) by the condition:
  \begin{displaymath}
    \tr p({}^1\!\cycle{s}{\bs}, {}^2\!\cycle{s}{\bs}, \ldots,  {}^n\!\cycle{s}{\bs})=0.
  \end{displaymath}
  However it is preferable to keep some geometrical meaning of
  constructed notions.

  An interesting observation is that in the matrix similarity of
  cycles~\eqref{eq:cycle-similarity} one may replace element
  \(g\in\SL\) by an arbitrary matrix corresponding to another cycle.
  More precisely the product
  \(\cycle{s}{\bs}\cycle[\tilde]{s}{\bs}\cycle{s}{\bs}\) is again the
  matrix of the form~\eqref{eq:FSCc-matrix} and thus may be associated
  to a cycle. This cycle may be considered as the reflection of
  \(\cycle[\tilde]{s}{\bs}\) in \(\cycle{s}{\bs}\).
  \begin{defn}
    \label{de:s-ortho}
    A cycle \(\cycle{s}{\bs}\) is s-orthogonal \emph{to} a cycle
    \(\cycle[\tilde]{s}{\bs}\) if the reflection of
    \(\cycle[\tilde]{s}{\bs}\) in \(\cycle{s}{\bs}\) is orthogonal
    (in the sense of Defn.~\ref{de:orthogonality}) to the real line.
    Analytically this is defined by:
    \begin{equation}
      \label{eq:s-orthog-def}
      \tr(\cycle{s}{\bs} \cycle[\tilde]{s}{\bs}\cycle{s}{\bs}\realline{s}{\bs})=0.
    \end{equation}
  \end{defn}
  Due to invariance of all components in the above definition
  s-orthogonality is a M\"obius invariant condition. Clearly
  this is not a symmetric relation: if \(\cycle{s}{\bs}\) is s-orthogonal to
  \(\cycle[\tilde]{s}{\bs}\) then \(\cycle[\tilde]{s}{\bs}\) is not
  necessarily s-orthogonal to
  \(\cycle{s}{\bs}\).
  
\begin{figure}[htbp]
  \includegraphics[]{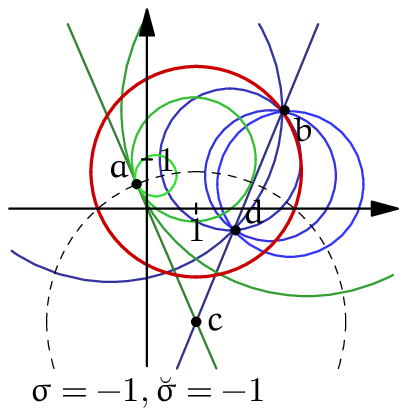}\hspace{1cm}
  \includegraphics[]{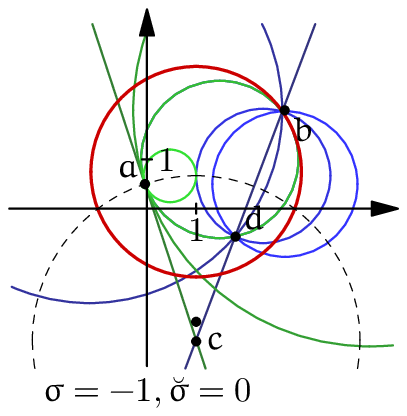}\hspace{1cm}
  \includegraphics[]{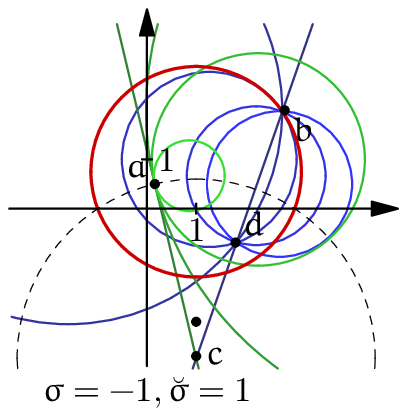}
  \caption[Orthogonality of the second kind]{Orthogonality of the
    second kind for circles. To highlight both
    similarities and distinctions with the ordinary orthogonality we
    use the same notations as that in Fig.~\ref{fig:orthogonality1}.}
  \label{fig:orthogonality2}
\end{figure}

  Fig.~\ref{fig:orthogonality2} illustrates s-orthogonality in the
  elliptic point space. By contrast with Fig.~\ref{fig:orthogonality1}
  it is not a local notion at the intersection points of cycles
  for all \(\bs\). However it may be again clarified in terms of the
  appropriate s-ghost cycle, cf. Thm.~\ref{th:ghost1}.
  \begin{thm}
    \label{th:ghost2}
    A cycle is s-orthogonal to a cycle \(C^{s}_{\bs}\) if its
    orthogonal in the traditional sense to its \emph{s-ghost cycle}
    \(\cycle[\tilde]{\bs}{\bs} = \cycle{\chi(\sigma)}{\bs}
    \Space[\bs]{R}{\bs} \cycle{\chi(\sigma)}{\bs}\), which is the
    reflection of the real line in \(\cycle{\chi(\sigma)}{\bs}\) and
    \(\chi\) is the \emph{Heaviside
      function}~\eqref{eq:heaviside-function}.  Moreover
    \begin{enumerate}
    \item \label{item:focal-centre-rel} \(\chi(\sigma)\)-Centre of
      \(\cycle[\tilde]{\bs}{\bs}\) coincides with the \(\bs\)-focus of
      \(\cycle{s}{\bs}\), consequently all lines s-orthogonal to
      \(\cycle{s}{\bs}\) are passing the respective focus.
    \item Cycles \(\cycle{s}{\bs}\) and \(\cycle[\tilde]{\bs}{\bs}\)
      have the same roots.
    \end{enumerate}
  \end{thm}
  Note the above intriguing interplay between cycle's centres
  and foci. Although s-orthogonality may look
  exotic it will naturally appear in the end of next Section again.

  Of course, it is possible to define another interesting higher order joint
  invariants of two or even more cycles.
  
  \section{Distance, length and perpendicularity}
  \label{sec:dist-lenght-perp}
  Geo\emph{metry} in the plain meaning of this word deals with \emph{distances}
  and \emph{lengths}. Can we obtain them from cycles?
  
\begin{figure}[htbp]
  \centering
  (a) \includegraphics[]{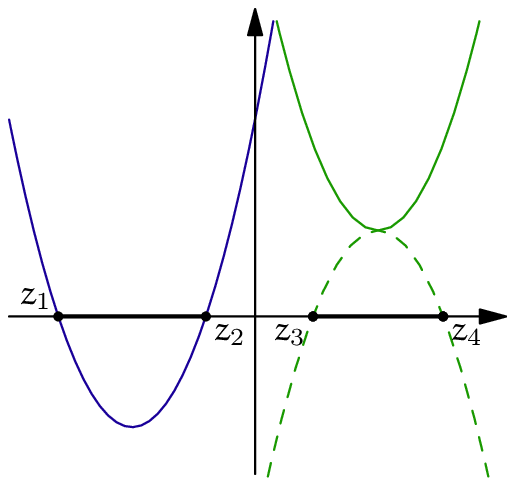}\hfill
  (b) \includegraphics[]{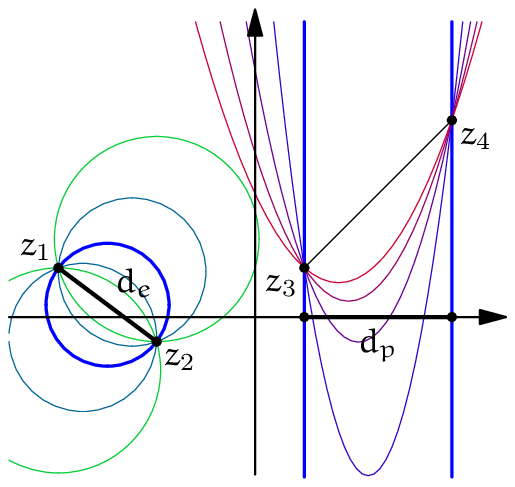}\hfill
  (c) \includegraphics[]{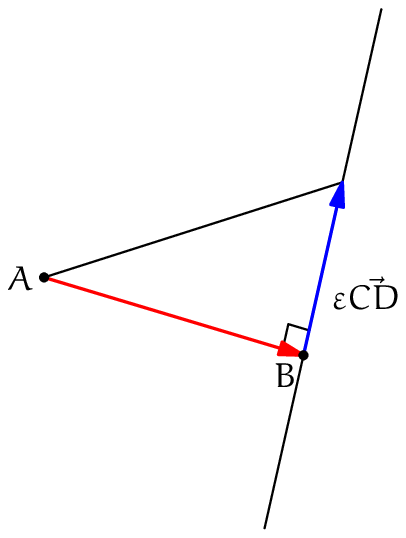}
  \caption[Radius and distance]{(a) The square of the parabolic
    diameter is the square of the distance between roots if they are
    real (\(z_1\) and \(z_2\)), otherwise the negative square of the
    distance between the adjoint roots (\(z_3\) and \(z_4\)).\\
    (b) Distance as extremum of diameters in elliptic (\(z_1\) and
    \(z_2\)) and parabolic (\(z_3\) and \(z_4\)) cases.\\
    (c) Perpendicular as the shortest route to a line.}
  \label{fig:distances}
\end{figure}

  We mentioned already that for circles normalised by the condition
  \(k=1\) the value \(\det
  \cycle{s}{\bs}=\scalar{\cycle{s}{\bs}}{\cycle{s}{\bs}}\) produces
  the square of the traditional circle radius. Thus we may keep it as the
  definition of the \emph{radius} for any cycle. But then we need to
  accept that in the parabolic case the radius is the (Euclidean) distance
  between (real) roots of the parabola, see
  Fig.~\ref{fig:distances}(a).

  Having radii of circles already defined we may use them for other
  measurements in several different ways. For example, the following
  variational definition may be used:

  \begin{defn}
    \label{de:distance}
    The \emph{distance} between two points is the extremum of
    diameters of all cycles passing through both points, see
    Fig.~\ref{fig:distances}(b).
  \end{defn}
  
  If \(\bs=\sigma\) this definition gives in all EPH cases the
  distance between endpoints of a vector \(z=u+\rmi v\) as follows:
  \begin{equation}
    \label{eq:eph-distance}
    d_{e,p,h}(u,v)^2=(u+\rmi v)(u-\rmi v)=u^2-\sigma  v^2.
  \end{equation}
  The parabolic distance \(d_p^2=u^2\), see
  Fig.~\ref{fig:distances}(b), algebraically sits between \(d_e\) and
  \(d_h\) according to the general principle~\eqref{eq:eph-class} and
  is widely accepted~\cite{Yaglom79}. However one may be unsatisfied
  by its degeneracy.

  An alternative measurement is motivated by the fact that a circle is
  the set of equidistant points from its centre. However the choice of
  ``centre'' is now rich: it may be either point from three
  centres~\eqref{eq:centres} or three foci~\eqref{eq:foci}.
  \begin{defn}
    \label{de:length}
    The \emph{length} of a directed interval \(\lvec{AB}\) is the radius
    of the cycle with its \emph{centre} (denoted by \(l_c(\lvec{AB})\))
    or \emph{focus} (denoted by \(l_f(\lvec{AB})\)) at the point \(A\)
    which passes through \(B\). 
  \end{defn}

  These definition is less common and have some unusual properties
  like non-symmetry: \(l_f(\lvec{AB})\neq l_f(\lvec{BA})\). However it
  comfortably fits the Erlangen program due to its
  \(\SL\)-\emph{conformal invariance}:

  \begin{thm}[\cite{Kisil05a}]
    Let \(l\) denote either the EPH distances~\eqref{eq:eph-distance}
    or any length from Defn.~\ref{de:length}. Then for
    fixed \(y\), \(y'\in\Space{R}{\sigma}\) the limit:
    \begin{displaymath}
      \lim_{t\rightarrow 0} \frac{l(g\cdot y, g\cdot(y+ty'))}{l(y,
        y+ty')}, \qquad
      \text{ where } g\in\SL, 
    \end{displaymath}
    exists and its value depends only from \(y\) and \(g\) and is
    independent from \(y'\).
  \end{thm}
  
  We may return from distances to angles recalling that in the
  Euclidean space a perpendicular provides the shortest root from a
  point to a line, see Fig.~\ref{fig:distances}(c). 
  \begin{defn}
    \label{de:perpendicular}
    Let \(l\) be a length or distance.  We say that a vector \(\lvec{AB}\) is
    \emph{\(l\)-perpendicular} to a vector \(\lvec{CD}\) if function
    \(l(\lvec{AB}+\varepsilon \lvec{CD})\) of a variable \(\varepsilon\) has a
    local extremum at \(\varepsilon=0\). 
  \end{defn}
  A pleasant surprise is that \(l_f\)-perpendicularity obtained
  thought the length from focus (Defn.~\ref{de:length}) coincides with
  already defined in Section~\ref{sec:higher-order-joint}
  s-orthogonality as follows from
  Thm.~\ref{th:ghost2}(\ref{item:focal-centre-rel}). It is
  also possible~\cite{Kisil08a} to make \(\SL\) action isometric in
  all three cases. 

  All these study are waiting to be  generalised to high dimensions
  and \wiki{Clifford_algebra}{Clifford algebras} provide a suitable
  language for this~\cite{Kisil05a,JParker07a}.

  \section{Erlangen program at large}
  \label{sec:erlangen-program-at}

  As we already mentioned the division of mathematics into areas is
  only apparent. Therefore it is unnatural to limit Erlangen program
  only to ``geometry''. We may continue to look for \(\SL\) invariant
  objects in other related fields. For example,
  transform~\eqref{eq:moebius} generates unitary
  representations on certain \(\FSpace{L}{2}\) spaces, cf.~\eqref{eq:moebius}:
  
  \begin{equation}
    \label{eq:hardy-repres}
    g^{-1}: f(x)\mapsto \frac{1}{(cx+d)^m}f\left(\frac{ax+b}{cx+d}\right).
  \end{equation}

  For \(m=1\), \(2\), \ldots the invariant subspaces of
  \(\FSpace{L}{2}\) are Hardy and (weighted) Bergman spaces of complex
  analytic functions.  All main objects of \emph{complex analysis}
  (Cauchy and Bergman integrals, Cauchy-Riemann and Laplace equations,
  Taylor series etc.) may be obtaining in terms of invariants of the
  \emph{discrete series} representations of
  \(\SL\)~\cite{Kisil02c}*{\S~3}.  Moreover two other series
  (\emph{principal} and \emph{complimentary}~\cite{Lang85}) play the
  similar r\^oles for hyperbolic and parabolic
  cases~\citelist{\cite{Kisil02c} \cite{Kisil05a}}.

  Moving further we may observe that transform~\eqref{eq:moebius} is
  defined also for an element \(x\) in any algebra \(\algebra{A}\)
  with a unit \(\mathbf{1}\) as soon as
  \((cx+d\mathbf{1})\in\algebra{A}\) has an inverse. If
  \(\algebra{A}\) is equipped with a topology, e.g. is a Banach
  algebra, then we may study a \emph{functional calculus} for element
  \(x\)~\cite{Kisil02a} in this way. It is defined as an intertwining
  operator between the representation~\eqref{eq:hardy-repres} in a
  space of analytic functions and a similar representation in a left
  \(\algebra{A}\)-module.

  In the spirit of Erlangen program such functional calculus is still
  a geometry, since it is dealing with invariant properties under
  a group action. However even for a simplest non-normal operator, e.g.
  a Jordan block of the length \(k\), the obtained space is not like a
  space of point but is rather a space of \(k\)-th
  \emph{jets}~\cite{Kisil02a}. Such non-point behaviour is oftenly
  attributed to \emph{non-commutative geometry} and Erlangen program
  provides an important input on this fashionable topic~\cite{Kisil02c}.

  Of course, there is no reasons to limit Erlangen program to \(\SL\)
  group only, other groups may be more suitable in different
  situations.  However \(\SL\) still possesses a big unexplored potential
  and is a good object to start with.

  \section{Acknowledgements}
  \label{sec:acknowledgements}

  I am in debt to the Editorial board of the \textsf{Notices of AMS} who
  spent enormous amount of time correcting and improving this paper.
  I am also grateful to Prof.~Troels Roussau Johansen who pointed out
  some more misprints and obscure places.

   Graphics for this article were created 
with the help of Open Source Software:
\small

\noindent MetaPost (\url{http://www.tug.org/metapost.html}),\\
Asymptote (\url{http://asymptote.sourceforge.net}),\\
and GiNaC (\url{http://www.ginac.de}).

\providecommand{\CPP}{\texttt{C++}} \providecommand{\NoWEB}{\texttt{noweb}}
  \providecommand{\MetaPost}{\texttt{Meta}\-\texttt{Post}}
  \providecommand{\GiNaC}{\textsf{GiNaC}}
  \providecommand{\pyGiNaC}{\textsf{pyGiNaC}}
  \providecommand{\Asymptote}{\texttt{Asymptote}} \newcommand{\noopsort}[1]{}
  \newcommand{\printfirst}[2]{#1} \newcommand{\singleletter}[1]{#1}
  \newcommand{\switchargs}[2]{#2#1} \newcommand{\irm}{\textup{I}}
  \newcommand{\iirm}{\textup{II}} \newcommand{\vrm}{\textup{V}}
  \providecommand{\cprime}{'} \providecommand{\eprint}[2]{\texttt{#2}}
  \providecommand{\arXiv}[1]{\href{http://arXiv.org/abs/#1}{arXiv:#1}}

\end{multicols}


\end{document}